\documentclass[10pt,a4paper]{article}
\usepackage{amssymb,amsmath,amsthm}

\usepackage{url,color}
\usepackage{times,lscape}
\def\de{{\rm d}}

\def\E{{\rm E}}
\def\Var{{\rm Var}}

\newtheorem{theorem}{Theorem}[section]

\newtheorem{remark}{Remark}[section]
\numberwithin{equation}{section}

\newtheorem{assumption}{Assumption}

\title{R\'{e}nyi information for ergodic diffusion processes}
\begin{document}

\author{Alessandro De Gregorio,  \,\,Stefano M. Iacus\footnote{\textbf{email:} alessandro.degregorio@unimi.it, stefano.iacus@unimi.it}\\
{\it Dipartimento di Scienze Economiche, Aziendali e Statistiche}\\
{\it    Via Conservatorio 7,  20122 Milan - Italy}
}

\maketitle

\begin{abstract}
In this paper we derive explicit formulas of the  R\'enyi information, Shannon entropy and Song measure for the invariant density of one dimensional ergodic diffusion processes. In particular, the diffusion models considered include the hyperbolic, the generalized inverse Gaussian,  the Pearson, the exponential familiy and a new class of skew-$t$ diffusions.
\end{abstract}



\section{Introduction}
For a continuous random variable $X$ with density $f$ the R\'enyi information of order $\alpha$ (R\'enyi 1959, 1961)
is defined as 
\begin{equation}\label{renyi}
\mathcal R_\alpha(f)=\frac{1}{1-\alpha}\log \int f^\alpha(x)\de x
\end{equation}
for $\alpha>0,\alpha\neq 1$ and 
\begin{equation}
\mathcal R_1(f) := \lim_{\alpha\to 1}\mathcal R_\alpha(f) = -\int f(x) \log f(x) \de x = - \E \log f(X)
\end{equation}
also known as Shannon or Kullback-Leibler entropy.
R\'enyi information is taken as a typical measure of complexity in the areas of physics, information theory  and engineering to describe \de ynamical or chaotic systems (see e.g Kurths {\it et al.}, 1995).
R\'enyi information, seen as a generalization of the Shannon entropy, is used to ``obtain different averaging of probabilities'' via the parameter $\alpha$ (see e.g. Song, 2001).
Considered as a function of $\alpha$, $\mathcal R_\alpha$ is also called the {\it spectrum} of the R\'enyi information and looking at its graphical plot is usually of interest. For some particular values of $\alpha$, it is used for specific tasks in different applications. For example, R\'enyi information of order $\alpha=2$  is used as a measure  of diversity in economics Hart (1975) and it is also the negative log of the well known limiting constant appearing in asymptotic efficiency lower bounds, i.e. $\int f^2(x) \de x$. 
Other applications to pattern recognition analysis (Vajda, 1968 and Ben Bassat, 1978), extreme value theory (see e.g. Falk {\it et al.}, 1994) and exponential models (Vajda and van der Meulen, 1998) have also appeared.
Recently Song (2001) put in evidence that not only the spectrum of the R\'enyi information is useful but also its gradient $\dot{\mathcal R}_\alpha = \frac{\partial}{\partial \alpha}\mathcal R_\alpha$. In particular, its value at $\alpha=1$ takes the following form
$$
\dot{\mathcal R}_1(f) = -\frac12  \Var (\log f(X))
$$
which means that the gradient of the spectrum evaluated at $\alpha=1$ is just the negative half of the variance of the loglikelihood whilst the $\mathcal R_1$ is its negative expected value.
Song proposed the following index
\begin{equation}\label{song}
\mathcal S(f) = -2\, \dot{\mathcal R}_1(f)\end{equation}
as an intrinsic measure of the shape of the density $f$. Notice that the information $\mathcal{R}_\alpha(f)$ is a measure of the dispersion of $f$ around its mean. The Song measure $\mathcal S$ is also a location and scale free functional in the sense that if $g(x)$ and $f(x) = g((x-\mu)/\sigma)/\sigma$ then $\mathcal S(f) = \mathcal S(g)$ and hence $\mathcal S$ it is a good measure of the shape of a distribution according to  Bickel and Lehman (1975). When $f$ has the fourth moment $\mu_4$, the Song measure provides similar information to the well known kurtosis measure $\beta_2 = \mu_4/\sigma^4$ but $\mathcal S$ measures more than what $\beta_2$ measures and it can be calculated also for distributions without $\mu_4$ (e.g. Student's $t$ with less than 4 degrees of freedom) or even without mean (e.g. the Cauchy distribution) and it is meaningful also for non symmetric distributions. Further, also in cases when two distribution have the same index $\beta_2$, the Song index may be different, e.g. the Student's $t$ with 6 degrees of freedom and the Laplace distribution have both $\beta_2=6$ but the Song measure is 1 for the Laplace distribution. It is approximatively 0.791 for the $t$ distribution. This indicates that their tails are not heavier in the same way. Then as proposed in Song (2001) the measure $\mathcal{S}(f)$ can be used to get a partial ordering respect to the tails of the distribution, i.e. given $f$ and $g$ two densities functions, the notation $f\prec g$ means that $\mathcal{S}(f) \leqslant\mathcal{S}(g)$. 

Liese and Vajda (1987) considered the following R\'enyi divergence between to densities $f$ and $g$
$$
D_\alpha(f,g) = \frac{1}{\alpha(1-\alpha)} \log \int \left(\frac{f(x)}{g(x)}\right)^\alpha g(x)\de x,\quad \alpha\neq 0,\alpha\neq 1,
$$
which can be used in hypotheses testing problems. The case $D_1(f,g)$ is the usual Kullback-Leiber  divergence. Suppose that $f(x,\theta)$ and $g(x,\theta)$ are two members of the same parametric family, $\theta\in \Theta$, and define the statistics $T_n^\alpha = 2 D_\alpha(f(\hat \theta_n), g(\theta_0))$ where $\hat\theta_n$ is a consistent estimator of $\theta$ and $\theta_0$ the value of $\theta$ under the null hypothesis, then  $T_n^\alpha$ can be used as test statistics and $T_n^1$ is just the usual likelihood ratio test statistics. Although a general theory is not yet available, several papers shows that for values of $\alpha>1$ the test statistics  $T_n^\alpha$ has often higher power than the corresponding likelihood ratio test statistics (see e.g. Rivas {\it et al.}, 2005, Morales {\it et al}, 1997, 2000, 2001).
Further, power divergences, as defined in Cressie and Read (1984), can be obtained from the R\'enyi divergence via the Box-Cox transformation
$$
\Psi_\alpha = \frac{ e^{\alpha(1-\alpha) D_\alpha} -1}{\alpha(1-\alpha)}\,.
$$
While likelihood inference for diffusion processes has now an extensive set of results for both continuous and discrete time observations, very few is known about the R\'enyi information and its related quantities for diffusion process. In particular, for continuos time observations from the diffusion process solution to the stochastic differential equation
$$
\de X_t^{(j)} = -a_j X_t^{(j)}\de t + \beta_t \de W_t, \quad X_0 \sim N(\mu_i, \sigma^2), \quad j=1,2,
$$
with $W_t$ a standard Wiener process and $a_j$ real constants, Vajda (1990) has obtained the explicit form of the R\'enyi divergence of the likelihood for all $\alpha>0$. In the discrete case Morales {\it et al.} (2005) obtained the explicit form of $T_n^\alpha$ for the model
$$
\de X_t = a_j \de t + b_j \de W_t, \quad X_0=x_0, \quad j=1,2\,,
$$
where $a_j$ and $b_j$ are real constants.
To the best of our knowledge, besides these two very simple cases, no other results are known for diffusion processes. In this work we consider homogeneous one dimensional diffusion processes solutions to the general stochastic differential equation
\begin{equation}
\de X_t = b(X_t) \de t + \sigma(X_t) \de W_t\,,
\label{eq:diff1}
\end{equation}
with some initial condition $X_0=x_0$ and state space
 $(l,r)$, with $-\infty \leq l \leq r \leq +\infty$. 
 
 In this paper we concentrate on the derivation of explicit formulas for the information measures $\mathcal R_\alpha(f)$, $\mathcal R_1(f)$ and $\mathcal S(f)$, where $f$ is the invariant density of several classes of ergodic diffusion processes including Pearson, hyperbolic, generalized inverse Gaussian, skew $t$ and exponential models. The laws considered are not included in other papers on R\'enyi informations for one dimensional densities like Song (2001), Zografos and Nadarajah (2003, 2005) which yet cover a large  variety of densities. At the same time, the densities considered here have been recently used to describe assets in financial markets because they describe well heavy tailed distributed returns (see e.g. Eberlein and Keller, 1995, Bibby and S\o rensen, 2001, Boyarchenko and Levendoroskii,2002). In the light of the above considerations, we think that the R\'enyi informations obtained here could be useful for researchers with particular attention to the financial applications to build new diagnostic tools.

\section{General considerations}
We introduce some basic assumptions on a diffusion process defined by \eqref{eq:diff1}. 
\begin{assumption}\label{ass:exun}
The drift and diffusion coefficient are such that the stochastic differential equation \eqref{eq:diff1} admits a unique weak solution $X_t$.
\end{assumption}
Let us introduce  the scale function and speed measure, defined respectively as
\begin{equation}
s(x) = \exp\left\{-2\int_{\tilde x}^x\frac{b(y)}{\sigma^2(y)}\de y\right\}\,,
\label{eq:scalem}
\end{equation}
with $\tilde x$ any value in the state space $(l,r)$,
and
\begin{equation}
m(x) = \frac{1}{\sigma^2(x)s(x)}\,.
\label{eq:speedm}
\end{equation}

\begin{assumption}\label{ass:ergodic}
Let $(l,r)$, with $-\infty \leq l \leq r \leq +\infty$ be the state space of the diffusion process $X_t$ solution to \eqref{eq:diff1} and assume that
$$\int_l^r m(x) \de x < \infty\,.
$$
Let, $x^*$ be  an arbitrary point  in the state space of $X$ such that either
$$
\int_{x^*}^r s(x) \de x = \int_l^{x^*} s(x) \de x = \infty\,.
$$
If one or both of the above integrals are finite, the corresponding boundary is assumed to be instantaneously  reflecting. 
\end{assumption}

If the Assumption \ref{ass:exun}-\ref{ass:ergodic} are satisfied, then exists a unique ergodic  process $X_t$ solution for the stochastic differential equation \eqref{eq:diff1}, with invariant law 
\begin{equation}
f(x) = \frac{m(x)}{G}\,,
\label{eq:invdens}
\end{equation}
where $G = \int m(x) \de x$. 

\begin{remark}\label{codomain}
The R\'{enyi} information defined by \eqref{renyi} for the invariant density function \eqref{eq:invdens} is equal to
\begin{equation}\label{eq:renyied}
\mathcal R_\alpha(f)=\frac{1}{1-\alpha}\left\{-\alpha\log G+\log\int_l^r\frac{1}{\sigma^{2\alpha}(x)}\exp{\left\{2\alpha\int_{\tilde{x}}^x\frac{b(y)}{\sigma^2(y)}\de y\right\}}\de x\right\}.
\end{equation}
We observe that $\lim_{\alpha\to\infty}\mathcal R_\alpha(f)=\log G$ and $\lim_{\alpha\to 0}\mathcal R_\alpha(f)=\log(\frac rl)$. Being the R\'{e}nyi information a psuedoconcave and monotonically decreasing function respect to $\alpha$ (see Bector and Bhatia, 1986), then we can conclude that $\inf_{\alpha>0}\mathcal R_\alpha(f)=\log G$ and $\sup_{\alpha>0}\mathcal R_\alpha(f)=\log(\frac rl)$. In other words the codomain of the function \eqref{eq:renyied} is $(\log G,\log(\frac rl)]$.
\end{remark}

\begin{remark}
The Shannon entropy for the distribution \eqref{eq:invdens} becomes
\begin{eqnarray*}
\mathcal{R}_1(f)&=&-\E\log f(x)\\
&=& \log G-	\E(\log \sigma^2(X_t))-  \E(\log s(X_t)),
\end{eqnarray*}
and for $b(x)=c\sigma^2(x), c\in \mathbb{R},$ we have that
$$\mathcal{R}_1(f)= \log G-E(\log \sigma^2(X_t)) +2c E(X_t-\tilde{X}).$$
\end{remark}

We point out that it is not always easy to calculate Shannon entropy directly as expected value of the logarithmic transformation of the density $f(x)$. So in this paper to obtain the esplicit form of $\mathcal{R}_1(f)$, we often use de l'H\^opital rule. Similar considerations can be done for the Song measure defined by \eqref{song}.

Furthermore we will denote by $\psi$ the digamma function, i.e. $\psi(x) = \Gamma'(x)/\Gamma(x)$ and its derivative by $\dot\psi(x) = \frac{\de}{\de x}\psi(x)$, while $B(a,b)=\frac{\Gamma(a)\Gamma(b)}{\Gamma(a+b)}$ is the Beta function.

\section{Pearson diffusions}
A Pearson diffusion (see Forman and Sorensen, 2006) is a stationary solution to a stochastic differential equation with mean reverting linear drift and squared diffusion coefficient of the form
\begin{eqnarray}\label{eq:pearson}
\de X_t&=&-\theta(X_t-\mu)\de t+\sigma(X_t)\de W_t\\
&=&-\theta(X_t-\mu)\de t+\sqrt{2\theta(aX_t^2+bX_t+c)}\de W_t\notag
\end{eqnarray}
where $\theta>0$ and the parameters $a,b$, and $c$ are such that the square root in \eqref{eq:pearson} is well defined. The scale function and the speed measure of $X_t$ are respectively
$$s(x)=\exp\left(\int_{x_0}^x\frac{y-\mu}{ay^2+by+c}dy\right),\quad m(x)=\frac{1}{2\theta s(x)(ax^2+bx+c)}$$
with $x_0$ a fixed point such that $ax_0^2+bx_0+c>0$. Let $(l,r)$ be the state space of the process such that $ax^2+bx+c>0$ for all $x\in (l,r)$. A unique ergodic solution to \eqref{eq:pearson} with initial value $x_0\in (l,r)$ exists if and only if the Assumption \ref{ass:ergodic} holds.  When the stationary process $X_t$ satisfying the stochastic differential equation \eqref{eq:pearson} exists, it is called Pearson diffusion because its invariant distribution belongs to the Pearson system. In fact we have that
$$\frac{\de m(x)}{\de x}=-\frac{(2a+1)x-\mu+b}{ax^2+bx+c}m(x)$$

We consider six particular cases as in Forman and Sorensen (2006) for \eqref{eq:pearson} which permit us to classify the ergodic Pearson diffusions by means of the squared diffusion coefficient.

\subsection{The Ornstein-Uhlenbeck diffusion: $\sigma(X_t)=\sqrt{2\theta}$}
For all $\mu \in \mathbb{R}$ exists a
unique solution of the equation \eqref{eq:pearson}. This process
is known as the Ornstein-Uhlenbeck diffusion and its invariant law $f$ is
a gaussian density with mean $\mu$ and variance 1. For $\alpha>0$, in this case the
R\'{e}nyi information, the Shannon entropy and the Song measure are respectively 
$$
\mathcal R_\alpha(f)=\frac12\left\{\log(2\pi)-\frac{\log\alpha}{1-\alpha}\right\},\quad \mathcal R_1(f)=\frac12\left\{1+\log(2\pi)\right\},\quad
\mathcal S (f)= 1\,.
$$
According to Remark \ref{codomain}, we note that in this case the R\'{enyi} information has codomain  $(\frac12\log(2\pi),\infty]$. 

\subsection{The  Cox-Ingersoll-Ross diffusion: $\sigma(X_t)=\sqrt{2\theta X_t}$}
We obtain the Cox-Ingersoll-Ross
process and the equation \eqref{eq:pearson} has a unique ergodic
solution on $(0,\infty)$ if and only if $\mu>1$. The invariant distribution is a Gamma
distribution with scale parameter 1 and shape parameter $\mu$.
Therefore for $\alpha>0$ (see Song, 2001)

$$
\begin{aligned}
\mathcal R_\alpha&=\frac{1}{1-\alpha}\left\{-\alpha\log\Gamma(\mu)-(\alpha(\mu-1)+1)\log\alpha+\log\Gamma(\alpha(\mu-1)+1)\right\},\\
\mathcal R_1&=\log\Gamma(\mu)-(\mu-1)\psi(\mu)+\mu, \\
\mathcal S &= \dot\psi(\mu ) (\mu -1)^2-\mu +2\,.
\end{aligned}
$$

\subsection{Pearson diffusion with $\sigma(X_t)=\sqrt{2\theta a(X_t^2+1)}$}
For all $a>0$ and $\mu\in \mathbb{R}$ the stochastic differential equation \eqref{eq:pearson} admits a unique ergodic solution. It's easy to verify that the scale function and speed measure are respectively 
$$s(x)=(1+x^2)^{\frac{1}{2a}}\exp\left\{-\frac{\mu}{a}\arctan x\right\}$$ 
and 
$$m(x)=\frac{1}{2\theta}(1+x^2)^{-\frac{1}{2a}-1}\exp\left\{\frac{\mu}{a}\arctan x\right\},$$ while the invariant density assumes the following form
\begin{equation}
f(x)=\frac{(1+x^2)^{-\frac{1}{2a}-1}\exp\left\{\frac{\mu}{a}\arctan x\right\}}{\int_{-\infty}^{+\infty} (1+x^2)^{-\frac{1}{2a}-1}\exp\left\{\frac{\mu}{a}\arctan x\right\}\de x}.
\end{equation}
If $\mu\neq 0$ the invariant law is skewed with tails decaying at the same rate as the 
$t$-distribution with $1+1/a$ degrees of freedom whilst for $\mu=0$ it is a scaled Student's $t$-distribution. This distribution is known as Pearson's type IV distribution, but a suitable name is the skew $t$-distribution with mean value $\mu$. 
The R\'{e}nyi information assumes the following form
\begin{eqnarray}\label{eq:renyip}
\mathcal{R}_{\alpha}(f)&=&\frac{1}{1-\alpha}\Bigg\{-\alpha\log\int_{-\infty}^{+\infty} (1+x^2)^{-\frac{1}{2a}-1}\exp\left(\frac{\mu}{a}\arctan x\right)\de x\\
&&+\log\int_{-\infty}^{+\infty} (1+x^2)^{-\frac{\alpha}{2a}-\alpha}\exp{\left(\alpha\frac{\mu}{a}\arctan x\right)}\de x \Bigg\}\notag\\
&=&\frac{1}{1-\alpha}\Bigg\{-\alpha\log\int_{-\pi/2}^{\pi/2}(\cos x)^{\frac{1}{a}}\exp\left(\frac{\mu}{a}x\right)\de x\notag\\
&&+\log\int_{-\pi/2}^{\pi/2}(\cos x)^{2\alpha\left(\frac{1}{2a}+1\right)-2}\exp\left(\frac{\alpha\mu}{a}x\right)\de x\Bigg\}\notag\\
&=&\frac{1}{1-\alpha}\Bigg\{-\alpha\log\int_{-\pi/2}^{\pi/2}(\cos x)^{\frac{1}{a}}\exp\left(-\frac{\mu}{a}x\right)\de x\notag\\
&&+\log\int_{-\pi/2}^{\pi/2}(\cos x)^{2\alpha\left(\frac{1}{2a}+1\right)-2}\exp\left(-\frac{\alpha\mu}{a}x\right)\de x\Bigg\}\notag
\end{eqnarray}
\begin{remark}
We observe that
$$\lim_{\alpha\to\infty}\mathcal{R}_{\alpha}(f)=\log\int_{-\pi/2}^{\pi/2}(\cos x)^{2\left(\frac{1}{2a}+1\right)-2}\exp\left(-\frac{\mu}{a}x\right)dx$$
and
$$\lim_{\alpha\to 0}\mathcal{R}_{\alpha}(f)=\log\pi$$
\end{remark}
Although the R\'{e}nyi information \eqref{eq:renyip} cannot be written in closed form, we are able to calculate it explicitly for some particular values of the parameter $a$ and of the order $\alpha$. 
We start with the next result.
\begin{theorem}\label{teopears}
Fixed $a>0$, for $\alpha=\frac{2a(m+1)}{1+2a}$ or $\alpha=\frac{2a(m+3/2)}{1+2a}$, with $m=1, 2,\ldots ,$ we obtain that
\begin{eqnarray}
\mathcal{R}_{\frac{2a(m+1)}{1+2a}}(f)&=&\frac{2a(m+1)}{2am-1}\log\int_{-\pi/2}^{\pi/2}(\cos x)^{\frac{1}{a}}\exp\left(-\frac{\mu}{a}x\right)\de x\\
&&+\frac{1+2a}{1-2am}\log\left(C(m,	\mu)\sinh\left(\frac{(m+1)\mu}{1+2a}\pi\right)\right)\notag
\end{eqnarray}
where $C(m,\mu)=\frac{2(2m)!}{\frac{2(m+1)\mu}{1+2a}\left(\left(\frac{2(m+1)\mu}{1+2a}\right)^2+2^2\right)\cdot\cdot\cdot \left(\left(\frac{2(m+1)\mu}{1+2a}\right)^2+(2m)^2\right)}$, and
\begin{eqnarray}\label{theo2}
\mathcal{R}_{\frac{2a(m+3/2)}{1+2a}}(f)&=&\frac{a(2m+3)}{a(2m+1)-1}\log\int_{-\pi/2}^{\pi/2}(\cos x)^{\frac{1}{a}}\exp\left(-\frac{\mu}{a}x\right)\de x\\
&&+\frac{1+2a}{1-a(2m+1)}\log\left(D(m,\mu)\cosh\left(\frac{(2m+3)\mu}{2(1+2a)}\pi\right)\right)\notag
\end{eqnarray}
where $D(m,\mu)=\frac{2(2m+1)!}{\left(\left(\frac{(2m+3)\mu}{1+2a}\right)^2+1^2\right)\cdot\cdot\cdot \left(\left(\frac{(2m+3)\mu}{1+2a}\right)^2+(2m+1)^2\right)}$.

\end{theorem}
\begin{proof}
For $\alpha=\frac{2a(m+1)}{1+2a}$, the espression \eqref{eq:renyip} gives  
\begin{eqnarray}\label{eq:renyip1}
\mathcal{R}_{	\frac{2a(m+1)}{1+2a}}(f)
&=&\frac{2a(m+1)}{2am-1}\log\int_{-\pi/2}^{\pi/2}(\cos x)^{\frac{1}{a}}\exp\left(-\frac{\mu}{a}x\right)\de x\\
&&+\frac{1+2a}{1-2am}\log\int_{-\pi/2}^{\pi/2}(\cos x)^{2m}\exp\left(-\frac{2(m+1)\mu}{1+2a}x\right)\de x\notag
\end{eqnarray}
By means of the formula 3.893(8) in Gradstheyn and Ryzhik (2007) 
\begin{eqnarray*}
&&\int_0^{\pi/2}(\cos x)^{2m}\exp{(-px)}\de x\\
&&=\frac{(2m)!}{p(p^2+2^2)\cdot\cdot\cdot (p^2+(2m)^2)}
\times\Biggl\{-e^{-p\pi/2}+1+\frac{p^2}{2!}+\frac{p^2(p^2+2^2)}{4!}+\cdots\\
&&\quad +\frac{p^2(p^2+2^2)\cdot\cdot\cdot(p^2+(2m-2)^2)}{2m!}\Biggr\}
\end{eqnarray*}
we derive
\begin{eqnarray}\label{eq:renyip2}
&&\int_{-\pi/2}^{\pi/2}(\cos x)^{2m}\exp\left(-\frac{2(m+1)\mu}{1+2a}x\right)\de x\\
&&=\int_{0}^{\pi/2}(\cos x)^{2m}\exp\left(-\frac{2(m+1)\mu}{1+2a}x\right)\de x\\
&&\quad+\int_{0}^{\pi/2}(\cos x)^{2m}\exp\left(\frac{2(m+1)\mu}{1+2a}x\right)\de x\notag\\
&&=\frac{(2m)! \left\{-e^{-\frac{(m+1)\mu}{1+2a}\pi}+e^{\frac{(m+1)\mu}{1+2a}\pi}\right\}}{\frac{2(m+1)\mu}{1+2a}\left(\left(\frac{2(m+1)\mu}{1+2a}\right)^2+2^2\right)\cdots\left(\left(\frac{2(m+1)\mu}{1+2a}\right)^2+(2m)^2\right)}\notag\\
&&=\frac{(2m)! 2\sinh\left(\frac{(m+1)\mu}{1+2a}\pi\right)}{\frac{2(m+1)\mu}{1+2a}\left(\left(\frac{2(m+1)\mu}{1+2a}\right)^2+2^2\right)\cdot\cdot\cdot \left(\left(\frac{2(m+1)\mu}{1+2a}\right)^2+(2m)^2\right)}\notag
\end{eqnarray}
Plugging in \eqref{eq:renyip2} inside the espression \eqref{eq:renyip1} the first part of the proof is concluded.
Analogously for $\alpha=\frac{2a(m+3/2)}{1+2a}$, we can show \eqref{theo2}. In fact
\begin{eqnarray*}
\mathcal{R}_{\frac{2a(m+3/2)}{1+2a}}(f)&=&\frac{a(2m+3)}{a(2m+1)-1}\log\int_{-\pi/2}^{\pi/2}(\cos x)^{\frac{1}{a}}\exp\left(-\frac{\mu}{a}x\right)\de x\\
&&+\frac{1+2a}{1-a(2m+1)}\log\int_{-\pi/2}^{\pi/2}(\cos x)^{2m+1}\exp\left(-\frac{(2m+3)\mu}{1+2a}x\right)\de x
\end{eqnarray*}
By means of the formula 3.893(10) in Gradstheyn and Ryzhik (2007)
\begin{eqnarray*}
&&\int_0^{\pi/2}(\cos x)^{2m+1}\exp{(-px)}\de x\\
&&=\frac{(2m+1)!}{(p^2+1)(p^2+3^2)\cdot\cdot\cdot (p^2+(2m+1)^2)}\\
&&\quad\times\left\{e^{-p\pi/2}+p \left\{1+\frac{p^2+1}{3!}+\cdots+\frac{(p^2+1)(p^2+3^2)\cdot\cdot\cdot (p^2+(2m-1)^2)}{(2m+1)!}\right\}\right\}
\end{eqnarray*}
we have that
\begin{eqnarray*}
&&\int_{-\pi/2}^{\pi/2}(\cos x)^{2m+1}\exp\left(-\frac{(2m+3)\mu}{1+2a}x\right)\de x\\
&&=\frac{(2m+1)! \left(e^{-\frac{(2m+3)\mu}{2(1+2a)}\pi}+e^{\frac{(2m+3)\mu}{2(1+2a)}\pi}\right)}{\left(\left(\frac{(2m+3)\mu}{1+2a}\right)^2+1^2\right)\cdots\left(\left(\frac{(2m+3)\mu}{1+2a}\right)^2+(2m+1)^2\right)}\\
&&=\frac{(2m+1)!2\cosh\left(\frac{(2m+3)\mu}{2(1+2a)}\pi\right)}{\left(\left(\frac{(2m+3)\mu}{1+2a}\right)^2+1^2\right)\cdot\cdot\cdot \left(\left(\frac{(2m+3)\mu}{1+2a}\right)^2+(2m+1)^2\right)}
\end{eqnarray*}
and \eqref{theo2} follows.

\end{proof}
The previous theorem permits us to derive the next result.
\begin{theorem}
Let $a=1$,  we have that
\begin{eqnarray}\label{eq:rpearson}
\mathcal{R}_{\frac{2m+2}{3}}(f)&=&\frac{2m+2}{2m-1}\log\left(\frac{2}{1+\mu^2}\cosh\left(\frac{\mu\pi} 2\right)\right)\\
&&+\frac{3}{1-2m}\log\left(C(m,	\mu)\sinh\left(\frac{(m+1)\mu}{3}\pi\right)\right)\notag
\end{eqnarray}
where $C(m,\mu)=\frac{2(2m)!}{\frac{(2m+2)\mu}{3}\left(\left(\frac{(2m+2)\mu}{3}\right)^2+2^2\right)\cdot\cdot\cdot \left(\left(\frac{(2m+2)\mu}{3}\right)^2+(2m)^2\right)}$, and
\begin{eqnarray}\label{eq:rpearson2}
\mathcal{R}_{\frac{2m+3}{3}}(f)&=&\frac{2m+3}{2m}\log\left(\frac{2}{1+\mu^2}\cosh\left(\frac{\mu\pi} 2\right)\right)\\
&&-\frac{3}{2m}\log\left(D(m,	\mu)\cosh\left(\frac{(2m+3)\mu}{6}\pi\right)\right)\notag
\end{eqnarray}
where $D(m,\mu)=\frac{2(2m+1)!}{\left(\left(\frac{(2m+3)\mu}{3}\right)^2+1^2\right)\cdot\cdot\cdot \left(\left(\frac{(2m+3)\mu}{3}\right)^2+(2m+1)^2\right)}$.
\end{theorem}
\begin{proof}
The quantities \eqref{eq:rpearson} and \eqref{eq:rpearson2} from Theorem \ref{teopears} by setting $a=1$  and noting that 
\begin{eqnarray}\label{eq:renyip3}
\int_{-\pi/2}^{\pi/2}\cos x\exp\left(-\mu x\right)\de x=\frac{2}{1+\mu^2}\cosh\left(\frac{\mu\pi} 2\right).
\end{eqnarray}

\end{proof}

\begin{theorem}
For $a=1$, the Shannon entropy is
\begin{equation}\label{eq:kl}
\mathcal{R}_1(f)=3\left\{\frac{\cosh\left(\frac{\mu\pi}{2}\right)}{1+\mu^2} -\frac32\log\Gamma'(2)-\frac12\mu\pi\tanh\left(\frac{\mu}{2}\pi\right)+\Gamma(2)\mu^2\right\}
\end{equation}
\end{theorem}
\begin{proof}
To calculate the Shannon entropy we consider the R\'{e}nyi information \eqref{eq:rpearson2} and observe that $\alpha\to 1$ if and only if $m\to 0$. Furthermore the relationship $\frac{\partial}{\partial\alpha}=\frac32\frac{\partial}{\partial m}$ holds. Then
\begin{eqnarray*}
\mathcal{R}_1(f)&=&\lim_{m\to 0} \mathcal{R}_{\frac{2m+3}{3}}(f)\\
&=&\lim_{m\to 0}-\frac94\Bigg\{-\frac23\log\int_{-\pi/2}^{\pi/2}\cos x\exp\left(-\mu x\right)\de x\\
&&+\frac{\frac{\partial}{\partial m}\left\{\int_{-\pi/2}^{\pi/2}(\cos x)^{2m+1}\exp\left(-\frac{(2m+3)\mu}{3}x\right)\de x\right\}}{\int_{-\pi/2}^{\pi/2}(\cos x)^{2m+1}\exp\left(-\frac{(2m+3)\mu}{3}x\right)\de x}\Bigg\}
\end{eqnarray*}
It's not hard to show that
\begin{eqnarray*}
&&\lim_{m\to 0}\frac{\partial}{\partial m}\left\{\frac{(2m+1)!2\cosh\left(\frac{(2m+3)\mu}{6}\pi\right)}{\left(\left(\frac{(2m+3)\mu}{3}\right)^2+1^2\right)\cdot\cdot\cdot \left(\left(\frac{(2m+3)\mu}{3}\right)^2+(2m+1)^2\right)}\right\}\\
&&=\frac{4\Gamma'(2)\cosh\left(\frac{\mu}{2}\pi\right)+\frac43\Gamma(2)\mu\pi\sinh\left(\frac{\mu}{2}\pi\right)-2\Gamma(2)\cosh\left(\frac{\mu}{2}\pi\right)\frac43\mu^2}{(1+\mu^2)},
\end{eqnarray*}
and by using again the result \eqref{eq:renyip3}, after some calculations we derive the expression \eqref{eq:kl}. 

\end{proof}

\subsection{Pearson diffusion with $\sigma(X_t)=\sqrt{2\theta a}X_t$}
If $a>0$ and $\mu>0$, we have an unique ergodic solution to the \eqref{eq:pearson}. In other words we obtain a diffusion process with invariant law
$$f(x)=\frac{a}{\mu\Gamma(1+\frac1a)}\left(\frac{\mu}{ax}\right)^{2+\frac1a}\exp{\left\{-\frac{\mu}{ax}\right\}},$$
that is an inverse Gamma distribution with shape $1+\frac1a$ and
scale parameter $\frac \mu a$. For all values of $\alpha$ such that $\alpha\left(2+\frac1a\right)-1>0$, it is not hard to prove that (see also Nadarajah and Zografos, 2003)

$$
\begin{aligned}
\mathcal R_\alpha(f)&=\log\left(\frac \mu a\right)+\frac{1}{1-\alpha}\left\{\left(1-	\left(2+\frac1a\right)\alpha\right)\log\alpha+\log\frac{\Gamma\left(\alpha\left(2+\frac1a\right)-1\right)}{\Gamma^\alpha \left(1+\frac 1a\right)}\right\}\\
\mathcal R_1(f)&=\log\left(\frac \mu a\right)+\log\Gamma\left(1+\frac1a\right)+\left(1+\frac1a\right)-\left(2+\frac1a\right)\psi\left(1+\frac1a\right)\\
\mathcal S (f)&=
	-\left(3 +	\frac1a\right) +\left(2+\frac1a\right)^2\dot\psi\left(1+\frac{1}{a}\right)
\end{aligned}
$$

\subsection{Pearson diffusion with $\sigma(X_t)=\sqrt{2\theta a X_t(X_t+1)}$}
If $a>0,\frac\mu a \geqslant 1$ we have that $m(x)=(
1+x)^{-\frac{\mu+1}{a}-1}x^{\frac\mu a-1}$ and the invariant density is a scaled $F$-distribution with $\frac{2\mu}{a}$ and $\frac2a+2$ degrees of freedom and scale parameter $\frac{\mu}{1+a}$, i.e.
$$f(x)=\frac{(
1+x)^{-\frac{\mu+1}{a}-1}x^{\frac\mu a-1}}{B\left(\frac\mu a, 1+\frac1a\right)}$$Hence
$$
\begin{aligned}
\mathcal R_\alpha(f)=&\frac{1}{1-\alpha}\biggl\{-\alpha\log B\left(\frac\mu a, 1+\frac1a\right)\\
&+\log B\left(\alpha\left(\frac\mu a-1\right)+1, \alpha\left(2+\frac1a \right)-1\right)\biggr\}\\
\mathcal R_1(f)=&\log B\left(\frac\mu a, 1+\frac1a\right)-\left(\frac\mu a-1\right)\psi\left(\frac\mu a\right)\\
&-\left(2+\frac1 a\right)\psi\left(1+\frac1 a\right)+\left(\frac\mu a+\frac1 a+1\right)\psi\left(\frac\mu a+\frac1 a+1\right)\\
\mathcal S(f)=&
\frac{\dot\psi\left(1+\frac{1}{a}\right) (2 a+1)^2+(a-\mu )^2 \dot\psi
\left(\frac{\mu }{a}\right)-(a+\mu +1)^2 \dot\psi\left(\frac{a+\mu +1}{a}\right)}{a^2}
\end{aligned}
$$

\subsection{Jacobi diffusion: $\sigma(X_t)=\sqrt{2\theta aX_t(X_t-1)}$}
If $a<0$,  and for all $\mu>0$ such that $\min(\mu,1-\mu)\geqslant -a$ exists an ergodic diffusion process on $(0,1)$ which satisfies the stochastic differential equation \eqref{eq:pearson} with speed measure $m(x)=(1-x)^{-\frac{1-\mu}{a}-1}x^{-\frac\mu a-1}$. The invariant law $f(x)$ is a Beta distribution with shape parameters $-\frac\mu a,-\frac{1-\mu}{a}$. This diffusion is also known as Jacobi diffusion. R\'enyi information measures become
$$
\begin{aligned}
\mathcal R_\alpha(f)=&\frac{1}{1-\alpha}\biggl\{-\alpha\log B\left(-\frac\mu a, -\frac{1-\mu}{a}\right)\\
&+\log B\left(-\alpha\left(\frac{\mu}{ a}+1\right)+1, -\alpha\left(\frac{1-\mu}{a}+1\right)+1\right)\biggr\}\\
\mathcal R_1(f)=&\log B\left(-\frac\mu a, -\frac{1-\mu}{a}\right)+\left(\frac\mu a+1\right)\psi\left(-\frac\mu a\right)\\
&+\left(\frac{1-\mu}{ a}+1\right)\psi\left(-\frac{1-\mu}{ a}\right)-\left(2+\frac1a\right)\psi\left(-\frac1a\right)\\
\mathcal S(f) =& 
-\dot\psi\left(-\frac{1}{a}\right) \left(2 +\frac{1}{a}\right)^2+\left(1-\frac{\mu -1}{a}\right)^2 \dot\psi
  \left(\frac{\mu -1}{a}\right)\\
  &+\left(1+\frac{\mu}{a} \right)^2 \dot\psi
   \left(-\frac{\mu }{a}\right)
\end{aligned}
$$

\section{Generalized inverse Gaussian diffusions}
We introduce (as in S\o rensen, 1997) the following stochastic differential equation
 \begin{equation}\label{sdegig}
 \de X_t=(\beta_1 X_t^{2\gamma-1}-\beta_2 X_t^{2\gamma}+\beta_3 X_t^{2(\gamma-1)})\de t+\lambda X_t^\gamma \de W_t,\quad X_0=x_0>0.
 \end{equation}
We assume that $\gamma\geqslant 0$, $\lambda>0$ and  set $\theta_i=2\beta_i/\lambda^2,i=1,2,3.$ The scale function is
\begin{equation*}
s(x)=C(x_0)x^{-\theta_1}\exp{\{\theta_2 x+\theta_3 x^{-1}\}}
\end{equation*}
while the speed measure becomes
\begin{equation*}
 m(x) =\frac{1}{C(x_0)}x^{\theta_1-2\gamma}\exp{\{-\theta_2 x-\theta_3 x^{-1}\}}.
 \end{equation*} 
 If the parameters satisfy the following conditions
 $$\theta_1\geqslant 1, \theta_2>0, \theta_3\geqslant 0$$
  $$1-2\gamma\leqslant\theta_1< 1, \theta_2>0, \theta_3> 0$$
   $$\theta_1< 1-2\gamma, \theta_2\geqslant0, \theta_3> 0$$
   the stochastic differential equation \eqref{sdegig} admits a unique ergodic weak solution (see S\o rensen, 1997) with invariant distribution 
  \begin{equation}\label{eq:invlawgig}
f(x)=\frac12 \left(\frac{\theta_2}{\theta_3}\right)^{\frac{\theta_1-2\gamma+1}{ 2}}\frac{x^{\theta_1-2\gamma}\exp\left\{-\theta_2x-\theta_3x^{-1})\right\}}{K_{\theta_1-2\gamma+1}\left(2\sqrt{\theta_2\theta_3}\right)}.
 \end{equation}
The law \eqref{eq:invlawgig} is a generalized inverse Gaussian density,
 where $K_\nu$ is the Bessel function of the third order wih index $\nu$ which has the following integral representation
 $$
 K_\nu(x)=\frac{1}{2}\int_0^\infty \exp\left\{-\frac x2( u+u^{-1})\right\}u^{-\nu-1}\de u,\quad \nu\in \mathbb{R}.$$
  \begin{theorem}
 The R\'{enyi} information for a generalized inverse Gaussian diffusion is
 \begin{equation}\label{eq:renyigig}
\mathcal{R}_{\alpha}(f)=- \log\left(\frac12\sqrt{\frac{\theta_2}{\theta_3}}\right)+\frac{1}{1-\alpha}\log \frac{K_{\alpha(\theta_1-2\gamma)+1}\left(2\alpha\sqrt{\theta_2\theta_3}\right)}{K_{\theta_1-2\gamma+1}^\alpha\left(2\sqrt{\theta_2\theta_3}\right)}.
\end{equation}
The Shannon entropy is equal to
\begin{eqnarray}\label{eq:shannongig}
\mathcal R_1 &=&  - \log\left(\frac12\sqrt{\frac{\theta_2}{\theta_3}}\right) + (\theta_1-2\gamma+1)+ \log K_{\theta_1-2
   \gamma +1}\left(2 \sqrt{\theta _2 \theta _3}\right)\\
&&-(\theta_1-2\gamma)
\frac{ \dot K_{\theta_1-2\gamma+1}\left( 2\sqrt{\theta_2\theta_3}\right)}{K_{ \theta _1-2 \gamma +1}\left(2 
   \sqrt{\theta _2 \theta _3}\right)} + 2\sqrt{\theta_2\theta_3}\frac{K_{\theta_1-2\gamma}\left(2\sqrt{\theta_2\theta_3}\right)}{K_{ \theta _1-2 \gamma +1}\left(2 
   \sqrt{\theta _2 \theta _3}\right)}\notag
\end{eqnarray}
with
$\dot K_\nu = \frac{\partial}{\partial \nu} K_\nu$. \\
 \end{theorem}
\begin{proof}
 To explict the R\'{e}nyi information we need to calculate the integral appearing in \eqref{renyi}. Hence, by means of the following integral representation for the Bessel function
  $$
 K_\nu(xz)=\frac{z^\nu}{2}\int_0^\infty \exp\left\{-\frac x2( u+z^2u^{-1})\right\}u^{-\nu-1}\de u,\quad \nu\in \mathbb{R},$$
and noting that $K_\nu(x)=K_{-\nu}(x)$,  we have for the invariant distribution \eqref{eq:invlawgig} the following result
 \begin{eqnarray}\label{eq:intgig}
 &&\int_0^\infty f^\alpha(x)\de x	\notag\\
 &&=\frac{1}{2^\alpha} \left(\frac{\theta_2}{\theta_3}\right)^{\frac{\alpha(\theta_1-2\gamma+1)}{ 2}}
 K^{-\alpha}_{\theta_1-2\gamma+1}\left(2\sqrt{\theta_2\theta_3}\right)\notag\\
 &&\quad\times \int_0^\infty x^{\alpha(\theta_1-2\gamma)}\exp\left\{-\frac\alpha2(2\theta_2x+2\theta_3x^{-1})\right\}\de x\notag\\
 &&=(y=2\theta_2 x)\notag\\
 &&= \left(\frac{\theta_2}{\theta_3}\right)^{\frac{\alpha(\theta_1-2\gamma+1)}{ 2}}\frac{
 K^{-\alpha}_{\theta_1-2\gamma+1}\left(2\sqrt{\theta_2\theta_3}\right)}{2^{\alpha(\theta_1-2\gamma+1)+1}}\notag\\
 &&\quad\times\int_0^\infty y^{\alpha(\theta_1-2\gamma)}\exp\left\{-\frac{\alpha}{2}(y+4\theta_2\theta_3y^{-1})\right\}\de y\notag\\
 &&=\left(\frac{\theta_2}{\theta_3}\right)^{\frac{\alpha(\theta_1-2\gamma+1)}{ 2}}\frac{
 K^{-\alpha}_{\theta_1-2\gamma+1}\left(2\sqrt{\theta_2\theta_3}\right)}{2^{\alpha-1}}\left(\frac{\theta_3}{\theta_2}\right)^{\frac{\alpha(\theta_1-2\gamma)+1}{2}}K_{\alpha(\theta_1-2\gamma)+1}\left(2\alpha\sqrt{\theta_2\theta_3}\right)\notag\\
 &&=\frac{1}{2^{\alpha-1}}\left(\frac{\theta_2}{\theta_3}\right)^{\frac{\alpha-1}{ 2}}
 \frac{K_{\alpha(\theta_1-2\gamma)+1}\left(2\alpha\sqrt{\theta_2\theta_3}\right)}{ K^{\alpha}_{\theta_1-2\gamma+1}\left(2\sqrt{\theta_2\theta_3}\right)}\notag
 \end{eqnarray}
Therefore we obtain the R\'{e}nyi information \eqref{eq:renyigig} by simple calculations. In order to obtain the Shannon entropy we use the following quantity
\begin{eqnarray*}
&&\frac{\partial}{\partial\alpha}A=\frac{\partial}{\partial\alpha}\left[
\log K_{\alpha  \left(\theta _1-2 \gamma \right)+1}\left(2 \alpha 
   \sqrt{\theta _2 \theta _3}\right) -\alpha  \log K_{\theta_1-2
   \gamma +1}\left(2 \sqrt{\theta _2 \theta _3}\right)
  \right] \\
  &&=\frac{\frac{\partial}{\partial\alpha}K_{\alpha  \left(\theta _1-2 \gamma \right)+1}\left(2 \alpha 
   \sqrt{\theta _2 \theta _3}\right)}{K_{\alpha  \left(\theta _1-2 \gamma \right)+1}\left(2 \alpha 
   \sqrt{\theta _2 \theta _3}\right)}-\log K_{\theta_1-2
   \gamma +1}\left(2 \sqrt{\theta _2 \theta _3}\right)
   \end{eqnarray*}
  
 Since $
\frac{\partial}{\partial x} K_\nu(x) =
-K_{\nu -1}(x)-\frac{\nu  K_{\nu }(x)}{x}
$ (see Gradshteyn and Ryzhik, 2007, formula 8.486(13)), it's easy to show that
\begin{eqnarray*}
&&\frac{\partial}{\partial\alpha}K_{\alpha  \left(\theta _1-2 \gamma \right)+1}\left(2 \alpha 
   \sqrt{\theta _2 \theta _3}\right)\\
   &&=(\theta-2\gamma)\dot K_{\alpha(\theta_1-2\gamma)+1}\left(2\alpha\sqrt{\theta_1\theta_2}\right)-2\sqrt{\theta_2\theta_3}K_{\alpha(\theta_1-2\gamma)}\left(2\alpha\sqrt{\theta_2\theta_3}\right)
   \\&&\quad-(\alpha(\theta_1-2\gamma)+1)K_{\alpha(\theta_1-2\gamma)+1}\left(2\alpha\sqrt{\theta_2\theta_3}\right).
   \end{eqnarray*}
   Therefore putting all together, by de l'H\^opital's rule, we have that 
$$
\begin{aligned}
\mathcal R_1(f) =&  - \log\left(\frac12\sqrt{\frac{\theta_2}{\theta_3}}\right) - \lim_{\alpha\to 1} \frac{\partial}{\partial \alpha} A\\
=& - \log\left(\frac12\sqrt{\frac{\theta_2}{\theta_3}}\right) + (\theta_1-2\gamma+1)+ \log K_{\theta_1-2
   \gamma +1}\left(2 \sqrt{\theta _2 \theta _3}\right)\\
&-(\theta_1-2\gamma)
\frac{ \dot K_{\theta_1-2\gamma+1}\left( 2\sqrt{\theta_2\theta_3}\right)}{K_{ \theta _1-2 \gamma +1}\left(2 
   \sqrt{\theta _2 \theta _3}\right)}\\
& + 2\sqrt{\theta_2\theta_3}\frac{K_{\theta_1-2\gamma}\left(2\sqrt{\theta_2\theta_3}\right)}{K_{ \theta _1-2 \gamma +1}\left(2 
   \sqrt{\theta _2 \theta _3}\right)}
\end{aligned}
$$

\end{proof}

The expression of the Song's measure can be obtained as well but it is too complicated and involves higher order derivatives of the Bessel's functions.
\begin{remark}
According to Remark \ref{codomain} the R\'{enyi} information for the generalized inverse Gaussian diffusions is a function with codomain $$\left(-\log\left(\frac12\sqrt{\frac{\theta_2}{\theta_3}}\right)+\log K_{\theta_1-2\gamma+1}\left(2\sqrt{\theta_2\theta_3}\right),\infty\right].$$
\end{remark}
\section{Hyperbolic diffusions}
Hyperbolic distributions have been extensively used in finance to model assets prices or log returns which deviates from the Gaussian law (see e.g. Eberlein and Keller, 1995 and Boyarchenko and Levendoroskii, 2002). Here we consider hyperbolic diffusion processes in the sense that their invariant law is an hyperbolic distribution.  The first simple version of hyperbolic process was introduced by Barndorff-Nielsen (1977) and further several generalizations and their applications in finance has been considered in Bibby and S\o rensen (1997). The solution of the following stochastic differential equation
\begin{equation}
\de X_t=\frac{\sigma^2}{2}\left\{\beta-\gamma\frac{X_t}{\sqrt{\delta^2+(X_t-\mu)^2}}\right\}\de t+\sigma
\de W_t
\end{equation}
is called hyperbolic diffusion (in the sense of Bibby and S\o rensen, 1997) with invariant law
\begin{equation}
f(x)=\frac{\sqrt{\gamma^2-\beta^2}}{2\gamma \delta
K_1(\delta\sqrt{\gamma^2-\beta^2})}\exp{\{-\gamma\sqrt{\delta^2+(x-\mu)^2}+\beta(x-\mu)\}},\quad
x\in\mathbb{R}.
\end{equation}

The next theorem represents the main result of this section. 
\begin{theorem}
For the hyperbolic diffusions we have that 
\begin{eqnarray}\label{eq:hyperbolicr}
\mathcal{R}_\alpha(f)&=&-\log\left(\frac{\sqrt{\gamma^2-\beta^2}}{2\gamma
\delta }\right)+\frac{1}{1-\alpha}\log\left(\frac{K_1(\alpha\delta\sqrt{\gamma^2-\beta^2})}{ K_1^\alpha(\delta\sqrt{\gamma^2-\beta^2})}\right)\\
\mathcal R_1(f)&=&
\log
K_1(\delta\sqrt{\gamma^2-\beta^2}) - \log\left(\frac{\sqrt{\gamma^2-\beta^2}}{2\gamma
\delta
}\right)+1\\
&&+ (\delta\sqrt{\gamma^2-\beta^2}) \frac{ K_0(\delta\sqrt{\gamma^2-\beta^2}) }{ K_1(\delta\sqrt{\gamma^2-\beta^2})},\notag
\end{eqnarray}
and
\begin{equation}
\begin{aligned}
\mathcal S(f)=&1 + \delta^2(\gamma^2-\beta^2) \Biggl\{ \frac{K_0\left( \delta \sqrt{\gamma ^2-\beta ^2} \right)K_2\left( \delta \sqrt{\gamma ^2-\beta ^2} \right)}{2K_1^2\left(\delta \sqrt{\gamma ^2-\beta ^2} \right)}\\
&-\frac{K_0^2\left(  \delta \sqrt{\gamma ^2-\beta ^2} \right)}{2K_1^2\left(\delta \sqrt{\gamma ^2-\beta ^2} \right)}\Biggr\}
\end{aligned}
\end{equation}
\end{theorem}
\begin{proof} We consider the integral

\begin{equation}
\mathcal I_\alpha=\left(\frac{\sqrt{\gamma^2-\beta^2}}{2\gamma \delta
K_1(\delta\sqrt{\gamma^2-\beta^2})}\right)^\alpha\int_{-\infty}^{+\infty}\exp{\{-\alpha\gamma\sqrt{\delta^2+(x-\mu)^2}+\alpha\beta(x-\mu)\}}\de x
\end{equation}
By means of a change of parametrization
$\varphi=\alpha(\gamma+\beta)$, $\psi=\alpha(\gamma-\beta)$, we
obtain that
\begin{eqnarray*}
&&\int_{-\infty}^{+\infty}
\exp{\{-\alpha\gamma\sqrt{\delta^2+(x-\mu)^2}+\alpha\beta(x-\mu)\}}\de x\\
&&=\int_{-\infty}^{+\infty}
\exp{\{-\alpha\gamma\sqrt{\delta^2+x^2}+\alpha\beta x\}}\de x\\
&&=\int_{-\infty}^{+\infty}
\exp{\left\{-\frac{1}{2}[\varphi(\sqrt{\delta^2+x^2}-x)+\psi(\sqrt{\delta^2+x^2}+x)]\right\}}\de x
\end{eqnarray*}
We note that the monotone transformation $\delta
y=\sqrt{\delta^2+x^2}+x$ maps $\mathbb{R}$ into $(0,\infty)$ with
inverse function $x=\frac{\delta}{2}\left(y-1/y\right)$. Therefore
\begin{eqnarray*}
&&\int_{-\infty}^{+\infty}
\exp{\left\{-\frac{1}{2}[\varphi(\sqrt{\delta^2+x^2}-x)+\psi(\sqrt{\delta^2+x^2}+x)]\right\}}\de x\\
&&=\frac{\delta}{2}\int_0^\infty
\left(1+\frac{1}{y^2}\right)\exp{\left\{-\frac\delta2\left(\frac\varphi
y+\psi y\right)\right\}}\de y\\
&&=\frac{\delta}{2}\left\{\frac1\psi\int_0^\infty\exp{\left\{-\frac\delta2\left(\frac{\varphi\psi}
{y}+ y\right)\right\}}\de y+\psi\int_0^\infty
y^{-2}\exp{\left\{-\frac\delta2\left(\frac{\varphi\psi} {y}+
y\right)\right\}}\de y\right\}\\
&&=\delta\left\{\frac{\sqrt{\varphi\psi}}{\psi}+\frac{\psi}{\sqrt{\varphi\psi}}\right\}K_1(\delta\sqrt{\varphi\psi})\\
&&=\frac{2\delta\gamma}{\sqrt{\gamma^2-\beta^2}}K_1\left(\alpha\delta\sqrt{\gamma^2-\beta^2}\right)
\end{eqnarray*}
and we can write 
\begin{equation}\label{eq:inithyperbolic}
\mathcal I_\alpha=\left(\frac{\sqrt{\gamma^2-\beta^2}}{2\gamma \delta
}\right)^{\alpha-1}\frac{K_1(\alpha\delta\sqrt{\gamma^2-\beta^2})}{K^\alpha_1(\delta\sqrt{\gamma^2-\beta^2})}.
\end{equation}
At this point by means of the integral $\mathcal{I}_\alpha$ the formula
\eqref{eq:hyperbolicr} easily follows.
To calculate the Shannon entropy we use de l'H\^opital rule to explicit the following limit
$$\lim_{\alpha\to1} \frac{1}{1-\alpha}\left\{\log K_1(\alpha\delta\sqrt{\gamma^2-\beta^2})-\alpha\log K_1(\delta\sqrt{\gamma^2-\beta^2})\right\}.$$
Hence by observing that
\begin{eqnarray*}
&&\frac{\partial}{\partial\alpha}
K_1(\alpha\delta\sqrt{\gamma^2-\beta^2})\\
&&=-\frac{\delta}{4}\int_0^\infty\left(\frac{\gamma^2-\beta^2}{x}+x\right)\exp{\left\{-\frac{\alpha\delta}{2}\left(\frac{\gamma^2-\beta^2}{x}+x\right)\right\}}\de x\\
&&=-\frac{\delta}{2}\sqrt{\gamma^2-\beta^2}\left\{K_0(\alpha\delta\sqrt{\gamma^2-\beta^2})+K_2(\alpha\delta\sqrt{\gamma^2-\beta^2})\right\}
\end{eqnarray*}
we get 
$$
\begin{aligned}
\mathcal R_1(f)=& - \log\left(\frac{\sqrt{\gamma^2-\beta^2}}{2\gamma
\delta
}\right)+\log
K_1(\delta\sqrt{\gamma^2-\beta^2})\\
&+\frac{\delta\sqrt{\gamma^2-\beta^2}\left\{K_0(\delta\sqrt{\gamma^2-\beta^2})+K_2(\delta\sqrt{\gamma^2-\beta^2})\right\}
}{2 K_1(\delta\sqrt{\gamma^2-\beta^2})}\\
=&
- \log\left(\frac{\sqrt{\gamma^2-\beta^2}}{2\gamma
\delta
}\right)+\log
K_1(\delta\sqrt{\gamma^2-\beta^2}) \\
&+1+ \delta\sqrt{\gamma^2-\beta^2} \frac{ K_0(\delta\sqrt{\gamma^2-\beta^2}) }{ K_1(\delta\sqrt{\gamma^2-\beta^2})}
\end{aligned}
$$
where in the last step we have used the relationship
$$
z \frac{K_0(z)+K_2(z)}{2K_1(z)} = 1+ z\frac{K_0(z)}{K_1(z)}.
$$

By the same calculations, one gets
$$
\begin{aligned}
\dot{\mathcal R}_\alpha(f)=&\frac{1}{(1-\alpha)^2}\Biggl\{
\frac{(\alpha -1) \delta \sqrt{\gamma ^2-\beta ^2}  K_0\left(\alpha  \delta
   \sqrt{\gamma ^2-\beta ^2} \right)}{K_1\left(\alpha  \delta
   \sqrt{\gamma ^2-\beta ^2} \right)}-\frac{1}{\alpha }+1\\
   &-
   \log \left( K_1\left( \delta\sqrt{\gamma ^2-\beta ^2} \right)\right)+\log
   \left(K_1\left(\alpha  \delta \sqrt{\gamma ^2-\beta ^2} \right)\right)
\Biggr\}
\end{aligned}
$$
which permits us to write
\begin{eqnarray*}
&&\mathcal S(f)\\
&&=\lim_{\alpha\to 1}\Bigg\{(\alpha-1)\delta \sqrt{\gamma ^2-\beta ^2}\frac{1}{K_1^2\left(\alpha  \delta \sqrt{\gamma ^2-\beta ^2} \right)}\Bigg[\frac{\partial}{\partial\alpha}K_0\left(\alpha  \delta \sqrt{\gamma ^2-\beta ^2} \right)K_1\left(\alpha  \delta \sqrt{\gamma ^2-\beta ^2} \right)\\
&&\quad-K_0\left(\alpha  \delta \sqrt{\gamma ^2-\beta ^2} \right)\frac{\partial}{\partial\alpha}K_1\left(\alpha  \delta \sqrt{\gamma ^2-\beta ^2} \right)\Bigg]+\frac{1}{\alpha^2}-\frac1\alpha\Bigg\}\frac{1}{-2(1-\alpha)}\\
&&=\delta^2(\gamma^2-\beta^2)\left(1-\frac{K_0^2\left(  \delta \sqrt{\gamma ^2-\beta ^2} \right)+K_0\left( \delta \sqrt{\gamma ^2-\beta ^2} \right)K_2\left( \delta \sqrt{\gamma ^2-\beta ^2} \right)}{2K_1^2\left(\delta \sqrt{\gamma ^2-\beta ^2} \right)}\right)+1
\end{eqnarray*}
\end{proof}
\begin{remark}
Since 
$$K_\nu(x)\simeq \sqrt{\frac{\pi}{2x}}e^{-x},\quad x\to\infty$$
and $K_\nu(0)=\infty$,
it's easy to see that
$$\lim_{\alpha\to\infty}\mathcal{R}_\alpha(f)=-\log\left(\frac{\sqrt{\gamma^2-\beta^2}}{2\gamma
\delta }\right)+\log K_1(\delta\sqrt{\gamma^2-\beta^2})$$
and
$$\lim_{\alpha\to 0}\mathcal{R}_\alpha(f)=\infty.$$
\end{remark}

\section{Skew Student's t-diffusions}
There are several skewed generalization of the Student's $t$-distribution which are widely used in financial context (see Aas and Haff, 2007, for a review). We introduce here a diffusion process which has the invariant distribution function that follows a skewed Student's $t$ distribution in the sense of Jones and Fad\de y (2003). We consider the process
solution of the following stochastic differential equation
\begin{equation}\label{eq:sdeskewt}
\de X_t=\frac{\sigma^2}{2}\frac{\gamma(\sqrt{\gamma+\beta+X^2_t}-X_t)-\beta(\sqrt{\gamma+\beta+X^2_t}+X_t)}{\gamma+\beta+X^2_t}\de t+\sigma
\de W_t
\end{equation}
with invariant law on the real line
\begin{equation}\label{denskew}
f(x)=\frac{\left(1+\frac{x}{\sqrt{\gamma+\beta+x^2}}\right)^{\gamma+\frac12}
\left(1-\frac{x}{\sqrt{\gamma+\beta+x^2}}\right)^{\beta+\frac12}}{B(\gamma,\beta)\sqrt{\gamma+\beta}2^{\gamma+\beta-1}},
\end{equation}
where $\gamma$ and $\beta$ are positive. We observe that when $\gamma=\beta$, $f(x)$ corresponds to a Student's t-distribution with $2\gamma$ degrees of freedom, whilst if $\gamma<\beta$ or $\gamma>\beta$ the density \eqref{denskew} is respectively positively or negatively skewed.  
\begin{remark}
To obtain the stochastic differential equation \eqref{eq:sdeskewt} we have borrowed the technique described in Bibby and S\o rensen (2003) for costructing diffusion processes with a specific invariant distribution. Let diffusion coefficient $\sigma(x)=\sigma$ be a real constant, the drift coefficient 
$$b(x)=\frac{\sigma}{2}\frac{\de}{\de x}	\log f(x),$$
where $f(x)$ is a density function, then \eqref{eq:sdeskewt} is an ergodic diffusion process with invariant law proportional to $f(x)$.
\end{remark}
\begin{theorem}
Under the conditions $\alpha(\gamma+\frac12)-\frac12>0$ and $\alpha(\beta+\frac12)-\frac12>0$ we obtain that

\begin{equation}\label{eq:renyiskew}
\mathcal R_\alpha (f)=\frac{1}{1-\alpha} \log\left(\frac{4^{\alpha-1}  B\left(\alpha\left(\gamma+\frac12\right)-\frac12,\alpha\left(\beta+\frac12\right)-\frac12\right)}{ \left(\sqrt{\beta+\gamma}\right)^{\alpha-1}B^\alpha(\gamma,\beta)}\right)
\end{equation}
and
\begin{eqnarray}\label{eq:shanskew}
\mathcal R_1(f)& =&-\log\left(\frac{4}{\sqrt{\gamma+\beta} B(\gamma,\beta)}\right)
-\left( \beta +\frac12\right) \psi(\beta )-\left( \gamma +\frac12\right) \psi
   (\gamma )\notag	\\
   &&+ (\beta +\gamma +1) \psi(\beta +\gamma
   )
   \end{eqnarray}
   \begin{equation}\label{eq:songskew}
\mathcal S(f) =
 \left( \beta +\frac12	\right)^2\dot\psi(\beta )+ \left( \gamma +\frac12	\right)^2\dot\psi(\gamma)-(\beta +\gamma +1)^2 \dot\psi (\beta +\gamma
   )
\end{equation}

\end{theorem}
\begin{proof} Since 
\begin{eqnarray*}
&&\int_{-\infty}^{+\infty}\left(1+\frac{x}{\sqrt{\gamma+\beta+x^2}}\right)^{\alpha(\gamma+\frac12)}
\left(1-\frac{x}{\sqrt{\gamma+\beta+x^2}}\right)^{\alpha(\beta+\frac12)}\de x\\
&&=(y=\frac{x}{\sqrt{\gamma+\beta+x^2}})\\
&&=\sqrt{\gamma+\beta}\int_{-1}^1(1+y)^{\alpha(\gamma+\frac12)-\frac32}(1-y)^{\alpha(\beta+\frac12)-\frac32}\de y\\
&&=(2w=1+y)\\
&&=2^{\alpha(\gamma+\beta+1)-2}\sqrt{\gamma+\beta}\int_{0}^1w^{\alpha(\gamma+\frac12)-\frac32}(1-w)^{\alpha(\beta+\frac12)-\frac32}\de w\\
&&=2^{\alpha(\gamma+\beta+1)-2}\sqrt{\gamma+\beta}\frac{\Gamma\left(\alpha(\gamma+\frac12)-\frac12\right)\Gamma\left(\alpha(\beta+\frac12)-\frac12\right)}{\Gamma\left(\alpha(\gamma+\beta+1)-1\right)}
\end{eqnarray*}
the R\'{e}nyi information \eqref{eq:renyiskew} immediately follows.
Starting by formula \eqref{eq:renyiskew} simple calculations permit us to derive the Shannon and Song measures \eqref{eq:shanskew}, \eqref{eq:songskew}.  
\end{proof}
\begin{remark}
We note that
$$\lim_{\alpha\to\infty}\mathcal R_\alpha(f)=\log B(\gamma,\beta) +\frac12\log(\gamma+\beta)-2\log2$$
and $\lim_{\alpha\to0}\mathcal R_\alpha(f)=\infty$.
\end{remark}
\section{Exponential families of diffusions}
Exponential families of diffusion processes have been introduced in S\o rensen (2007) and are solutions to the following stochastic differential
equation
\begin{equation}\label{eq:exp}
\de X_t = \sum_{i=1}^p \beta_i b_i(X_t)\de t + \lambda \sigma(X_t)\de W_t\,,
\end{equation}
with $\lambda >0$ and $v>0$ and $b_i$ functions such that a unique weak solution exists.
These models are not exponential families of stochastic processes in the sense of  K\"uchler and S\o rensen (1997).
By using the reparametrization $\theta_i = \beta_i/\lambda^2$ and the functions
$$
T_i(x) = 2\int_{x_0}^x \frac{b_i(y)}{\sigma^2(y)}\de y
$$
for some $x_0$ in the state space of $X$, the scale measure has the following representation
$$
s(x,\theta) = \exp\left\{-\sum_{i=1}^p \theta_i T_i(x) \right\}
$$
and the speed measure is
$$
m(x,\theta) = \sigma(x)^{-2}  \exp\left\{\sum_{i=1}^p \theta_i T_i(x) \right\}\,.
$$
If we denote by $\phi(\theta) = \int m(x,\theta) \de x$ the normalizing constant,
the invariant distribution of this class of models takes the canonical form of an exponential family, i.e.
$$
f(x) = \exp\left\{
\sum_{i=1}^p \theta_i T_i(x) - 2\log \sigma(x) - \phi(\theta)
\right\}\,.
$$
These models are ergodic (and hence the invariant density $f(x)$ exists) under proper conditions on the parameter  $\theta_i$, $i=1, \ldots, p$. In particular, it is need that $\theta\in \Theta \cap \Theta_1$, where
$$\Theta=\left\{\theta\in \mathbb{R}^p:\int_l^rm(x,\theta)\de x<\infty	\right\}$$
and
$$\Theta_1=\left\{\theta\in \mathbb{R}^p:\int_l s(x,\theta)\de x=\int^r s(x,\theta)\de x=\infty	\right\}.$$
For a general exponential diffusion the R\'{e}nyi information is equal to
\begin{equation}
\mathcal{R}_\alpha(f)=\frac{1}{1-\alpha}\left\{-\alpha\phi(\theta)+\log\int \frac{1}{\sigma^{2\alpha}(x)}\exp{\left\{\alpha\sum_{i=1}^p \theta_i T_i(x)\right\}}\de x\right\}.
\end{equation}
\begin{remark}
Although it is not possible to obtain explicitly  $\mathcal R_\alpha$, $\mathcal R_1$ and $\mathcal S$ in the general case,  it is not hard to see that exponential families of diffusions include, as particular cases, the Cox-Ingersoll-Ross, the Generalized Inverse Gaussian, the hyperbolic and the Ornstein-Uhlenbeck diffusions (see e.g. S\o rensen, 2007).
\end{remark}

\section*{Conclusions}
In this paper we have been able to obtain explicit forms of $\mathcal R_\alpha$, $\mathcal R_1$ and $\mathcal S$  for several classes of diffusion processes widely used in many fields including finance. 
These results, may be useful as a reference to develop new diagnostic tools and estimation procedures. The family of exponential diffusion processes, introduced in the last section, includes most of the other models for which we obtained explicit results as particular cases, so this family may be used as a building block to define goodness of fit tests for embedded families of diffusion processes based on the R\`enyi divergences. This will be the topic of further investigations.

 \end{document}